\documentstyle[12pt]{amsart}
\newcommand{\ck}{{\cal K}}

\newcommand{\cad}{{\cal D}}
\newcommand{\cf}{{\cal F}}
\newcommand{\ce}{{\cal E}}
\newcommand{\cg}{{\cal G}}
\newcommand{\cll}{{\cal L}}

\newcommand{\ds}{\displaystyle}
\newtheorem{th}{Theorem}
\newtheorem{lem}{Lemma}
\newtheorem{df}{Definition}
\newtheorem{rem}{Remark}

\makeatletter
\@addtoreset{equation}{section}

\makeatother

\begin{document}

\title
[SOME FUNCTION SPACES RELATED TO THE BROWNIAN MOTION ON FRACTALS]
{SOME FUNCTION SPACES RELATED TO THE BROWNIAN MOTION ON SIMPLE NESTED FRACTALS}
\author{Katarzyna Pietruska-Pa\l uba}
\address{\hskip-\parindent
  Katarzyna Pietruska-Pa\l uba\\
  Institute of Mathematics\\
  University of Warsaw\\
  ul. Banacha 2\\
  02-097 Warsaw, Poland}
\email{kpp\@mimuw.edu.pl}

\thanks{Research at MSRI
was supported in part by NSF grant DMS-9701755.}

\begin{abstract}
 In this paper we identify the domain of the Dirichlet form
associated with the Brownian  motion on simple nested fractals
with an integral Lipschitz space. This result generalizes
such an identification on the Sierpi\'nski gasket,
carried on by Jonsson in \cite{Jon}.
\end{abstract}

\maketitle

\section{Introduction}
Since the mid-eighties there has been an outburst of papers concerned
with the Brownian motion on fractal spaces, dealing with
the existence/uniqueness problem as well as investigating properties
of the resulting process(es). A good source of references for those
papers is \cite{Bar1}.

The first set to be investigated has been the Sierpi\'nski gasket,
by far the simplest representant of the `simple nested fractals' class
(see  \cite{Bar-Per},\cite{Gol},\cite{Kus}). The first constructions
 were purely 
probabilistic and the Brownian motion was obtained as a limit
(in distribution) of appropriately scaled random walks on lattices,
approximating
the gasket (see  \cite{Bar-Per},  for example). Later it became clear that
simpler approach to the problem is the one that uses Dirichlet forms rather
than random walks (see \cite{Fuk-Shi}). Both constructions
were then generalized to the class of `simple nested fractals' (see Sec. 2.1 
for definition) --- probabilistically in \cite{Kum}, \cite{Lin} and from
the Dirichlet point of view in \cite{Fuk}.

 The Dirichlet form associated with
the semigroup it generates
 is a local regular Dirichlet form on $L^2$ 
on the fractal with respect to the Hausdorff measure, which is invariant
under those isometries.  As Brownian motion on a simple
nested fractal is unique (up to a trivial time-rescaling), both approaches
are equivalent.

In the early paper  of Barlow and Perkins
\cite{Bar-Per} it was proven that  all functions
in the domain of the generator of the Brownian motion
have good continuity properties (they are H\"older continuous) and that 
they cannot be extended beyond their fractal domain as differentiable
functions. 

Jonsson in  \cite{Jon} gave a precise description of the domain of the
Dirichlet
form associated with the Brownian motion on the Sierpi\'nski gasket;
it turns out that this domain coincides
 with certain integral Lipschitz space on the fractal in question.
Those spaces as well as closely related  Besov spaces on
general sets  are
analyzed in depth in \cite{Jon-Wal}.

Jonsson conjectured that analogous result should hold for a much larger class
of fractals. This is true indeed --- by refinement of his methods 
we were able to extend the characterisation of the Dirichlet form domain
to the class of all simple nested fractals.
The advantage of simple nested fractals is that the Dirichlet form
can be obtained as a limit of finite-dimensional forms (see Sec. \ref{dir}),
which is the main tool in our proof.

\section{Preliminaries. Simple nested fractals and Dirichlet forms}
Throughout the paper, $\#K$ denotes the  cardinality of the set $K$
and $C(K)$ --- the space of continuous real-valued functions on $K.$
Moreover, generic
constants whose values are irrelevant to our purposes are denoted by $c$
and they usually change from line to line.

\subsection{Simple nested fractals}\label{snf}
The class of simple nested fractals was first 
introduced by T. Lindstr\o m
in \cite{Lin}. Our short exposition will be mostly
based on that paper.

Let $N\geq 1$ be fixed. A transformation $\psi:{\bf R}^N\to
{\bf R}^N$ is called a {\it similitude with scaling factor} $L$
($L>1)$ 
if $\psi(x)=\frac{1}{L}U(x)+v,$ where $U$ is an isometry of ${\bf R}^N$
and $v\in{\bf R}^N$ is fixed.

Suppose that $\psi_1,...,\psi_M, M\geq 2,$ are given similitudes with
common scaling factor $L.$ In view of Hutchinson's result
(see \cite{Fal}, \cite{Hut}, \cite{Lin}) there exists a unique 
nonempty compact set $\ck\subset {\bf R}^N$ such that
\[\ck=\bigcup_{i=1}^M\psi_i(\ck).\]
It is called the \underline{\em  self-similar fractal} generated by the
family of similitudes $\psi_1,...\psi_M.$

There are several explicit
methods of constructing $\ck.$ Let us describe one of them, based on
 the set of fixed points of the similitudes $\psi_i.$
 
 As
all the similitudes $\psi_i$ are contractions,  
there exist $x_1,...,x_M$ such that $\psi_i(x_i)=x_i$ (these
fixed points are not necessarily distinct). Denote by $F$ the
collection of those fixed points.

\begin{df}
$x\in F$ is called an  \underline{essential fixed point}
 of the system $\psi_1,...
\psi_M$ if there exists another fixed point $y\!\in\! F$
 and two different
transformations $\psi_i,\psi_j$
 such that $\psi_i(x)=\psi_j(y).$
\end{df}
Informally speaking, the essential fixed points are points through which 
parts of the fractal meet.
Denote by $V_0$ the set of all essential fixed points.

\vskip 1mm
\noindent\underline{\em Condition 1.} $\# V_0\geq 2.$
\vskip 3mm

\noindent{\bf Example 1.} {\em The Sierpi\'nski gasket in ${\bf R}^2.$}
 Three
transformations:
$\psi_1(x)=\frac{x}{2},$
$\psi_2(x=\frac{x}{2}+(\frac{1}{2},0),$
$\psi_3(x)=\frac{x}{2}+(\frac{1}{4},\frac{\sqrt{3}}{4}).$ The system $(\psi_1,
\psi_2,\psi_3)$ has three fixed points: $(0,0), (1,0), (\frac{1}{2},
\frac{\sqrt{3}}{2})$ and all of them are essential (see fig. 1).

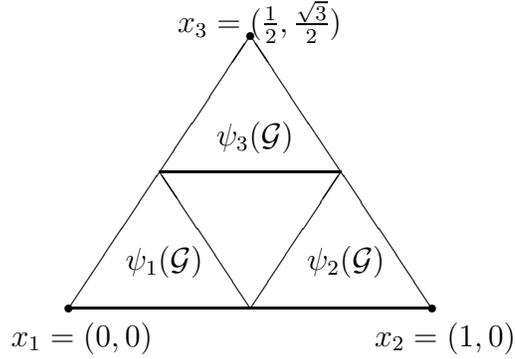
\begin{figure}
\addtolength\unitlength{-0.1mm}
\begin{picture}(300,200)(-130,0)
\put(40,30){\makebox(8,4)[lb]{$\psi_1({\cal G})$}}
\put(136,30){\makebox(8,4)[lb]{$\psi_2({\cal G})$}}
\put(88,96){\makebox(8,4)[lb]{$\psi_3({\cal G})$}}
\put(10,10){\line(1,0){192}}
\put(10,10){\circle*{4}}
\put(202,10){\circle*{4}}
\put(106,154){\circle*{4}}
\put(12,-10){\makebox(10,4)[b]{$x_1=(0,0)$}}
\put(204,-10){\makebox(10,4)[b]{$x_2=(1,0)$}}
\put(106,156){\makebox(10,4)[b]{$x_3=(\frac{1}{2},\frac{\sqrt 3}{2})$}}
\put(58,82){\line(1,0){96}}
\put(10,10){\line(2,3){96}}
\put(106,10){\line(2,3){48}}
\put(202,10){\line(-2,3){96}}
\put(106,10){\line(-2,3){48}}
\end{picture}
\caption{Transformations that define the Sierpi\'nski gasket ${\cal G}$}
\end{figure}

\noindent{\bf Example 2.} {\em The Lindstr\o m snowflake} is an invariant
set for the system of seven transformations $\psi_1,...\psi_7. $ Of the
fixed points $x_1,...,x_7,$ all are essential fixed points, except for
$x_7$ (see fig. 2)

\begin{figure}
\begin{picture}(300,200)(-65,0)
\put(70,10){\line(2,0){40}}
\put(150,10){\line(2,0){40}}
\put(30,70){\line(2,0){200}}
\put(30,130){\line(2,0){200}}
\put(70,190){\line(2,0){40}}
\put(150,190){\line(2,0){40}}
\put(10,100){\line(2,3){20}}
\put(50,160){\line(2,3){20}}
\put(50,40){\line(2,3){100}}
\put(110,10){\line(2,3){100}}
\put(190,10){\line(2,3){20}}
\put(230,70){\line(2,3){20}}
\put(70,10){\line(-2,3){20}}
\put(30,70){\line(-2,3){20}}
\put(150,10){\line(-2,3){100}}
\put(210,40){\line(-2,3){100}}
\put(250,100){\line(-2,3){20}}
\put(210,160){\line(-2,3){20}}
\put(70,10){\circle*{3}}
\put(65,2){\makebox(8,4){$x_1$}}
\put(190,10){\circle*{3}}
\put(195,2){\makebox(8,4){$x_2$}}
\put(250,100){\circle*{3}}
\put(258,96){\makebox(8,4){$x_3$}}
\put(190,190){\circle*{3}}
\put(195,192){\makebox(8,4){$x_4$}}
\put(70,190){\circle*{3}}
\put(58,192){\makebox{$x_5$}}
\put(10,100){\circle*{3}}
\put(-2,96){\makebox(8,4)[l]{$x_6$}}
\put(130,100){\circle*{2}}
\put(126,92){\makebox(8,4){$x_7$}}
\put(78,43){\makebox(8,4)[lt]{$\psi_1({\cal H})$}}
\put(158,43){\makebox(8,4)[lt]{$\psi_2({\cal H})$}}
\put(78,163){\makebox(8,4)[lt]{$\psi_5({\cal H})$}}
\put(158,163){\makebox(8,4)[lt]{$\psi_4({\cal H})$}}
\put(38,103){\makebox(8,4)[lt]{$\psi_6({\cal H})$}}
\put(198,103){\makebox(8,4)[lt]{$\psi_3({\cal H})$}}
\put(118,113){\makebox(8,4)[lt]{$\psi_7({\cal H})$}}
\end{picture}
\caption{Transformations that define the Lindstr\o m snowflake $\cal H$}
\end{figure}
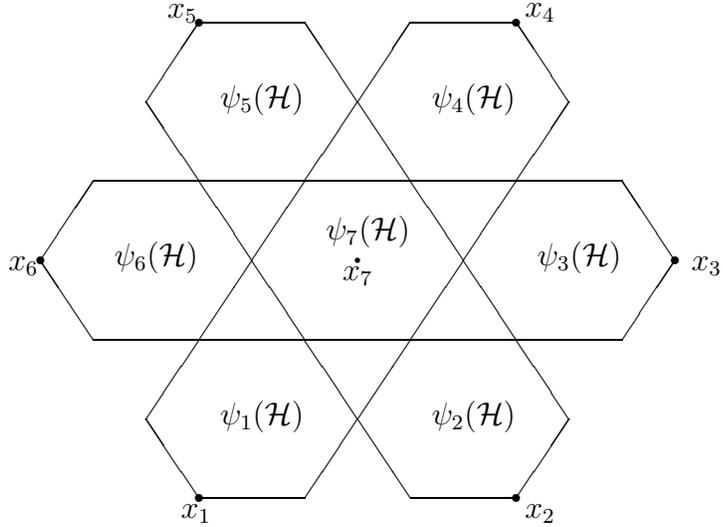

\noindent For $A\subset {\bf R}^N$  let
$\Psi(A)=\bigcup _{i=1}^{M} \psi_i(A)$ and $\Psi^n=
 \underbrace{\Psi\circ\cdots\circ\Psi}_{n\, \mbox{\scriptsize  times }}.$
Then we define `all the vertices at the $n$-th level' as
\[V_n =\Psi^n(V_0),\;\; n=1,2...\] 

 \noindent $(V_n)_n$ is an increasing
family of sets (it follows from $V_0\subset\Psi(V_0)).$ Let
\[V_{\infty}=\bigcup_{n=1}^{\infty}V_n.\]
The following holds true (see \cite{Fal}, p. 120):

\begin{th}
 $\ck=\overline{V_{\infty}},$ the closure taken in the Euclidean
topology of ${\bf R}^N.$
\end{th}

\noindent Next, we impose the following 
\vskip 1mm
\noindent
{\em \underline{Condition 2} (Open Set Condition).}
The family $\psi_1,...\psi_M$ satisfies the {\it\underline{open set}
\underline{condition}}
if there exists an open, bounded, nonempty set $U\subset {\bf R}^N$ such
that $\psi_i(U)\cap\psi_j(U)=\emptyset$ for $i\neq j$ and
\[\Psi(U)=\bigcup_{i=1}^{M} \psi_i(U)\subset U.\]

\vskip 1mm

If the open set condition is satisfied, then the Hausdorff dimension 
of $\ck$ can be easily calculated  (see Th. 8.6 of \cite{Fal})
and is equal to
\begin{equation}\label{i1}
d_f(\ck)=\frac{\log M}{\log L}\leq N
\end{equation}

For instance, the Hausdorff dimension of the two-dimensional Sierpi\'nski
gasket equals to $\frac{\log 3}{\log 2}$ and 
of the Lindstr\o m snowflake --- to $\frac{\log 7}{\log 3}.$ 
 $M$ is sometimes 
 called the mass-scaling factor, and $L$ -- the length-scaling factor
 for the fractal.

\begin{df}
 Let $m\geq 1$ be given.   A set
of the form $\psi_{i_1}\circ\dots\circ\psi_{i_m}(\ck)$ is called an
\underline{$m$-symplex}
($m$-symplices are just scaled down copies of $\ck).$ Collection
of all $m$-symplices will be denoted by $\cf_m.$  A set of the form 
$\psi_{i_1}\circ\dots\circ\psi_{i_m}(V_0)$ is called an 
\underline{$m$-cell}.
For an $m$-symplex $S=\psi_{i_1}\circ\dots\circ\psi_{i_m}(\ck),$
let $V(S)=\psi_{i_1}\circ\dots\circ\psi_{i_m}(V_0)$ be the set of
its vertices.
\end{df}

\noindent 
$(i_1,...,i_m)$ is called the \underline{\it address} of the symplex
$S=\psi_{i_1}\circ...\circ\psi_{i_m}(\ck).$
Similarly, if $x=\psi_{i_1}\circ\dots\circ\psi_{i_m}(x_0),$ $x_0\in V_0,$
then we call $(i_1,...i_m)$  the \underline{\it address} of the point
$x.$
Observe that addresses of symplices are uniquely determined, addresses
of points from $V_{\infty}$ in general are not.

\begin{df}
Let $m\geq 1 $ be fixed.
$x,y\in V_{\infty}$ are called \underline{$m$-neighbours} if there exists an
$m$-symplex
$S$ such that $x,y\in V(S)$ (in particular $x,y\in V_m).$
\end{df}

\noindent
After these definitions we can formulate the remaining three conditions.

\vskip 2mm
\noindent{\it \underline{Condition 3.} Nesting.} For each $m\geq 1,$ and
$S,T\in \cf_m,$ $S\cap T=V(S)\cap V(T).$ This condition is in
fact equivalent to 
\[\forall _{i,j\in\{1,...,M\}}\;\; \psi_i(\ck)\cap\psi_j(\ck)=\psi_i(V_0)
\cap\psi_j(V_0).\]
\vskip 1mm
\noindent{\it \underline{Condition 4.} Connectivity.} Define the graph
 structure 
$E_{(1)}$ on $V_1$ as follows:
{ we say that $(x,y)\in E_{(1)}$ if
$x$ and $y$ are $1$-neighbours.}
 Then we require the graph $(V_1,E_{(1)})$
to be connected.
\vskip 1mm
\noindent{\it \underline{Condition 5.} Symmetry.} For $x,y\in V_0$ let
$R_{x,y}$
be the reflection in the hyperplane bisecting the segment $[x,y].$
Then we stipulate that
\[
\forall _{i\in\{1,...,M\}}\forall _{x,y\in V_0,\;x\neq y}
\exists _{j\in\{1,...,M\}}\;\;  R_{x,y}(\psi_i(V_0))=\psi_j(V_0)
\]
(natural reflections map 1-cells  onto 1-cells).

\begin{df} Fractals satisfying conditions 1-5 are called \underline{%
simple nested} \underline{ fractals}.
\end{df}

\begin{rem}{\em The most restrictive of those assumptions is nesting, as
this condition requires the parts of the fractal to meet only through
vertices. Sierpi\'nski carpet, for example, does not fall into this category.
However, this condition is crucial for the construction of the Dirichlet form.}
\end{rem}

\begin{lem} 1.
 \begin{equation}\label{i2}
 \exists _{c_1,c_2 >0}\;\;\mbox{\it such that} \;\;\forall_{ m\geq 1}\;\;
 c_2 M^m \leq \#V_m\leq c_1 M^m.
\end{equation}
2.\begin{equation}\label{i3}
\forall_{ m\geq 1}\forall _{S\in\cf_m }\forall _{n\geq m}\;\;\#[V_n\cap S]=
\#
V_{n-m}.
\end{equation}
\end{lem}

\noindent{\bf Proof.}  Denote by $v_m$ the cardinality of $V_m.$ 
Obviously $v_0<v_1,$ and also\linebreak $v_1 <M\cdot v_0$ 
(as we deal with essential fixed points
only).

Let $k_0= Mv_0-v_1\geq 1.$ Then for each $m$ we have (nesting)
\[v_{m+1}=M\cdot v_m -k_0,\] so that
\[v_m = M^m\left(v_0-\frac{k_0}{M-1}\right)+\frac{k_0}{M-1}.\]
In view of definition of $k_0$ we have $(M-1)v_0-k_0=v_1-v_0 >0$
and (\ref{i2}) follows.
(\ref{i3}) is obvious.
$\Box$

\vskip 2mm
\noindent Finally, for $m\geq 1$
let $\mu_m$ be the normalized counting measure
on $V_m,$ i.e.
\begin{equation}\label{i4}
\mu_m(\cdot)=\frac{1}{\# V_m}\sum_{x\in V_m}\delta_{\{x\}}(\cdot).
\end{equation}

It is not hard to see that
measures $\mu_m$ converge weakly towards the $d_f$-dimensional
Hausdorff measure restricted to the fractal $\ck,$ $\mu.$

\subsection{Dirichlet forms on simple nested fractals}\label{dir}
The material of this section is based on  \cite{Bar} and \cite{Fuk}.

Suppose that the nested fractal $\ck,$ generated by the family
of similitudes $\psi_1,...,\psi_M$ is given. 
Let $A=(a_{x,y})_{x,y\in V_0}$ be a \underline{\it conductivity matrix}
on $V_0$, i.e. a real-valued matrix that satisfies
\begin{enumerate}
\item $\forall_{ x\neq y}\;\;a_{x,y}\geq 0\;\mbox{ and } a_{x,y}=a_{y,x},$
\item $\forall _x\;\sum_{y} a_{x,y}=0$ (so that $a_{x,x}\leq 0)$.
\end{enumerate}

We assume that $A$ is \underline{\it irreducible} i.e. the graph
$(V_0, E_{(A)})$ is connected, where $E_{(A)}$ is the graph structure on
$V_0$
related to the matrix $A:$
\[ \mbox{for } x,y\in V_0 \;\;
(x,y)\in E_{(A)} \Leftrightarrow a_{x,y}>0.\]

\noindent\underline{\it Dirichlet form on $V_0$ associated  with $A$}
 is defined as follows:
for $f\in C(V_0)$ (i.e. for any $f:\!V_0\to {\bf R},$ as $V_0$ is finite)
\begin{equation}\label{i5}
\ce_A^{(0)}(f,f)=\ce_A(f,f)=\frac{1}{2}\sum_{x,y\in
V_0}a_{x,y}(f(x)-f(y))^2\geq 0.
\end{equation}

\noindent We introduce two operators acting on Dirichlet
forms, related to the geometrical structure of $\ck.$

\noindent 1. {\it \underline{Reproduction}.} For $f\in C(V_1)$ let
\[\widetilde{\ce}^{(1)}_A(f,f)=\sum_{i=1}^M 
\ce^{(0)}_A(f\circ\psi_i,f\circ\psi_i).\]
The mapping $\ce^{(0)}_A\mapsto\widetilde{\ce^{(1)}_A}$ is called the {\it
\underline{reproduction
 map}}
and will be denoted by $\cal R.$
\vskip 2mm
\noindent 2. {\it \underline{Decimation}.} Given a symmetric form $\ce$ on
$V_1,$
define its restriction to $V_0,$ $\ce|_{V_0},$ by:\\
for $f:\!V_0\to {\bf R},$ 
\begin{equation}\label{i6}
\ce|_{V_0}=\inf\{\ce(g,g):\! g:V_1\to{\bf R}\mbox{ and } g|_{V_0}=f\}.
\end{equation}
This mapping is called the {\it\underline{decimation map}}
 and is denoted by ${\cal D}e.$

Denote by $\bf G$ the symmetry group of $V_0,$ i.e. the group of 
transformations generated by $\{R_{x,y}: x,y\in V_0\}.$

The following holds true (Lindstr\o m \cite{Lin} for the existence part,
 Sabot \cite{Sab} for the uniqueness part):

\begin{th}\label{thNDHS}
There exist a unique number $\rho = \rho(\ck)$ and a unique (up to
a multiplicative constant) irreducible conductivity matrix $A $ on $V_0,$
invariant under $\bf G$ and such that
\[({\cal D}e\circ {\cal R})(\ce_A)=\frac{1}{\rho}(\ce_A).\]
$A$ is called the symmetric nondegenerate harmonic structure on $V_0$
(symmetric NDHS, or just NDHS). It follows then that $a_{x,y}>0$ for
$x\neq y$ and $\rho>1.$
\end{th}

This nondegenerate harmonic structure is of particular interest for defining
the Dirichlet form on the complete fractal.

Suppose from now that the conductivity matrix $A$ is equal to the NDHS
from theorem \ref{thNDHS} which allows us to drop the subscript `$A$'.

For $f\in C(V_m)$ define $\widetilde{\ce}^{(m)}(f,f)$ 
in a manner similar to  (\ref{i6}), i.e.
\[\widetilde{\ce}^{(m)}(f,f)=\sum_{i_1,...,i_m} \ce^{(0)}(f\circ\psi_{i_1}\circ
...\circ\psi_{i_m}, f\circ\psi_{i_1}\circ...\circ\psi_{i_m})\]
 and further
\[\ce^{(1)}(f,f)\stackrel{def}{=} \rho \widetilde{\ce}^{(1)}(f,f),\;\;
\ce^{(m)}(f,f)\stackrel{def}{=} \rho^m \widetilde{\ce}^{(m)}(f,f).\]
With this notation, as a consequence of nesting
 we have:

\vskip 1mm
\begin{th}\label{thform}
\[ \forall _{f\in L^2(\ck,\mu)}\forall _{m=0,1,2,...}
\;\;\ce^{(m+1)}(f,f)\geq \ce^{(m)}(f,f).\]
Set
\[\cad=\cad(\ce)=\{f\in L^2(\ck,\mu):\sup_m\ce^{(m)}(f,f)<+\infty\}\]
and for $f\in\cad$
\begin{equation}\label{i7}
\ce(f,f)=\lim_{m\to\infty}\ce^{(m)}(f,f).
\end{equation}
Then $\cad\subset C(\ck).$ Equation (\ref{i7}) defines a Dirichlet form 
on $L^2(\ck,\mu)$ (a regular, symmetric, closed local  bilinear form)
that is invariant under ${\bf G}.$
\end{th}
In particular `$\ce$ is closed' means that $\cad$ is closed under the
norm
\begin{equation}\label{i8}
\|f\|_{\ce}\stackrel{def}{=} \sqrt{\ce(f,f)+\|f\|_2^2}.
\end{equation}
We can write the forms $\ce^{(m)}$ in a more convenient way:
\begin{equation}\label{aa}
{\ce}^{(m)}(f,f)=\rho^m\sum_{S\in\cf_ m} \widetilde{\ce}^{(m)}_S(f,f),
\end{equation}
where for $\cf_m\ni S=\psi_{i_1}\circ...\circ\psi_{i_m}(\ck)$
\[\begin{array}{lr}
\widetilde{\ce}^{(m)}_S(f,f) & =
\ce^{(0)}(f\circ\psi_{i_1}\circ\dots\circ\psi_{i_m},
f\circ\psi_{i_1}\circ\dots\circ\psi_{i_m})=\\[2mm]
{} & =\frac{1}{2}\ds\sum_{x,y\in V(S)}a_{x,y}^{(m)}(f(x)-f(y))^2;
\end{array}\]
$a_{x,y}^{(m)}$ being the `$m$-th level conductivities' between $x,y\in
V(S):$
if $x,y$ are $m$-neighbours and
$x=\psi_{i_1}\circ\dots\circ\psi_{i_m}(\bar  x),$ 
$y=\psi_{i_1}\circ\dots\circ\psi_{i_m}(\bar y),$ then
\[a^{(m)}_{x,y}=a_{\bar x, \bar y}\] 
and if $x,y$ are not $m$-neighbours we set $a^{(m)}_{x,y}=0.$

\subsubsection{Significance of the constant $\rho.$}\label{constants}
By analogy with the electrical circuit theory, $\rho$ is called
the \underline{\em resistance scaling factor}.
Denote 
\[d_w=d_w(\ck)\stackrel{def}{=}\frac{\log (M\rho)}{\log L}>1.\]
 As the constant $\rho$
is uniquely determined, this number can be considered 
to be one of characteristic numbers of the fractal $\ck$
and is called  `the \underline{\it walk dimension} of $\ck$'.
One has
\begin{equation}\label{i9}
\rho=\frac{L^{d_w}}{L^{d_f}}.
\end{equation}
Next,
 \[d_s=d_s(\ck)\stackrel{def}{=} \frac{2d_f}{d_w}\]
 is the \underline{\em spectral
dimension} of the set $\ck.$ In view of $\rho >1$ we have $d_s <2,$
which will be decisive for our results in the sequel.
For example, for the  $n$-dimensional
Sierpi\'nski gasket $\cal G$ the walk dimension equals to
$\frac{\log(n+3)}{\log 2}$ and its spectral dimension --- to
$\frac{2\log 3}{\log 5}<2.$

\subsection{Lipschitz spaces on $d$-sets}
(see Jonsson-Wallin \cite{Jon-Wal}, Jonsson \cite{Jon})

Let $F$ be an arbitrary closed subset of ${\bf R}^N $ ($N\geq 1$ is fixed)
and $d$ a real number, $0<d\leq N.$ $B(x,r)$ denotes the closed ball with
centre $x$ and radius $r.$ 
A positive Borel measure $\mu$ with support $F$   is called 
a \underline{\it $d$-measure} on $F$ if there exist constants $c_1$ and $c_2$ 
such that
\[
\forall_{x\in F}\forall_{ 0<r\leq 1}\;\;
{c_1 r^{ d}\leq \mu(B(x,r))\leq c_2 r^{d}}.
\]
If $F$ is a support or a $ d$-measure, then $F$ is called a {\em 
\underline{$ d$-set}.}
Nested fractal $\ck$ with Hausdorff dimension $d_f$ is a $d_f$-set
and the restriction of the $d_f$-dimensional Hausdorff measure to
$\ck$ is a $d_f$-measure (it follows from \cite{Fal}, pp. 120-121)

If $F$ is a closed subset of ${\bf R}^N$ (not necessarily
 a $d$-set),
 and $\gamma >0$ is a real number, then 
$Lip (\gamma, F)$ is the set of those bounded
functions $f:\!F\to {\bf R}$
that satisfy the H\"{o}lder condition with exponent $\gamma:$
\begin{equation}\label{eqi0}
\exists _{M>0}
\forall _{x,y\in F, |x-y|\leq 1}\;\;
 |f(x)|\leq M,\;\;
 |f(x)-f(y)|\leq M|x-y|^{\gamma}.
\end{equation}
The norm of $f$ in $Lip(\gamma, F)$ is the infimum of all the
possible constants $M$ in (\ref{eqi0}).

If, for some $d>0,$
 $F\subset {\bf R}^N$ is a $d$-set and $\mu$ is a
  $d$-measure on $F,$
 then one defines the following `integral' Lipschitz
 spaces on $F:$\\
for given $\alpha >0, \;\;1\leq p,q\leq +\infty$
 the Lipschitz space $Lip(\alpha,p,q)(F)$ 
 is the collection of all those functions $f\in L^p(F,\mu)$
 for which, when $p,q <\infty$
 \[\|(a_m^{(p)}(f))\|_{l_q}=\left(\sum_{m=1}^{\infty}(a^{(p)}_m(f))^q
\right)^{\frac{1}{q}}<\infty,\]
where coefficients $a_m^{(p)}(f)$ are defined as
\begin{equation}\label{i10}
a_m^{(p)}(f)=2^{m\alpha}\left(2^{md}\int\!\int_{|x-y|<\frac{c_0}{2^m}}
|f(x)-f(y)|^pd\mu(x)d\mu(y)\right)^{\frac{1}{p}}
\end{equation}
for some $c_0>0.$  For $p=\infty$ or $q=\infty$ 
  take the $L_{\infty}$
or $l^{\infty}$ norms, as needed.

The space $Lip(\alpha,p,q)(F)$ is
a Banach space with the norm
\begin{equation}\label{i11}
\|f\|_{Lip}=\|f\|_{L^p(\ck,\mu)}+\|(a_m^{(p)}(f))\|_{l_q}.
\end{equation}

Note that if we replace the constant $c_0$ in (\ref{i10})
 with another one 
 we obtain the same function space
and the norm changes to an equivalent one. The same remains true if
we replace `2' by some other number greater than 1, see lemma \ref{lem1}
below.

Finally let us formulate the following property (Corollary 2, p. 498
of \cite{Jon}):
\begin{th}\label{continuity}
For $0<d\leq n,$ $1\leq p,q\leq\infty$ and $\alpha$ such that
$\alpha-\frac{d}{p}=\gamma \in(0,1),$ $Lip(\alpha,p,q)(F) $ is continuously
imbedded in $Lip (\alpha-\frac{d}{p},F).$ In particular
every $f\in Lip(\alpha,p,q)(F)$ is H\"{o}lder continuous
with exponent $\gamma.$
\end{th}

Jonsson (\cite{Jon}) proves that on the $N$-dimensional Sierpi\'nski 
gasket $\cg\subset{\bf R}^N$ ($\cg$ is a $d_f$-set with
 $d_f=\frac{\log (N+3)}{\log 2}$) the domain of the Dirichlet form defined by 
(\ref{i7}) coincides with the Lipschitz space
\[
Lip(\alpha,2,\infty)(\cg),\;\;\mbox{with}\;\alpha=\frac{d_w}{2}=
\frac{\log(N+3)}{2\log 2}
\]
 and the two norms:
the Lipschitz norm and the Dirichlet norm are equivalent. He conjectured 
that a similar result should hold for a larger class of fractals.
In this paper we prove that this conjecture 
holds  for all  simple nested fractals.

\section{The main theorem}

The natural scaling factor for the fractal $\ck$ is $L=L(\ck),$
and in general $L\neq 2.$ Therefore it is convenient to substitute
$L$ for 2 in the definition of  Lipschitz coefficients, (\ref{i10}). This
minor
alteration changes the Lipschitz
norm to an equivalent one. This is true in general: if $F\subset {\bf R}^N$
is a $ d-$set, $\alpha >0, 0<p< \infty$    ($\alpha$ is considered fixed) let
\begin{equation}\label{eq1}
b_m^{(p)}(f)= L^{m\alpha}\left(L^{m d}\int\!\int_{|x-y|<\frac{c_0}{L^m}}
|f(x)-f(y)|^p d\mu(x)d\mu(y)\right)^{\frac{1}{p}}
\end{equation}
(similarly for $p=\infty).$

\noindent We have:
\begin{lem}\label{lem1}
Let $F\subset {\bf R}^N$ be a $d$-set, $d\!\in(0,N],$ $\mu $ --- a $d$-measure
on $F.$ Then
there exists a constant $D>0$ such that for each $f\in L^p(F, \mu)$
and for each $0<q\leq\infty$
\[\frac{1}{D}\| (a_m^{(p)}(f))\|_{l_q}\leq \| (b_m^{(p)}(f))\|_{l_q}
\leq D \| (a_m^{(p)}(f))\|_{l_q}.
\] 
\end{lem}

\noindent{\bf Proof.} For $m\geq 1,$ let $n(m)$ be the unique integer
satisfying
\[2^{n(m)}\leq L^m < 2^{n(m)+1}\]
(for large $m,$ $n(m)$ is the integer part of $m \frac{\log 2}{\log
L})$
so that
\[\begin{array}{ll}
b_m^{(p)}(f) & \leq (2^{n(m)+1})^{\alpha}\left((2^{n(m)+1})^{ d}
\ds\int\!\ds\int _{|x-y|<\frac{c_0}{2^{n(m)}}}|f(x)-f(y)|^p
d\mu(x)d\mu(y)\right)^{\frac{1}{p}}\\
{} & = 2^{\alpha+\frac{ d}{p}}a_{n(m)}^{(p)}(f)
\end{array}\]
and
\[\begin{array}{rl}
\| (b_m^{(p)}(f))\|_{l_q}^q=\sum_m(b_m^{(p)}(f))^q & \leq (2^{\alpha+
\frac{ d}{p}})^q \sum_m (a_{n(m)}^{(p)}(f))^q \\[2mm]
{} & \leq (2^{\alpha+\frac{ d}{p}})^q\|(a_m^{(p)}(f))\|_{l_q}^q
\end{array}\]
(similarly for $q=\infty).$ The opposite inequality follows by symmetry.
$\Box$

\vskip 2mm
From now on we restrict our attention to $F=\ck,$ $ d= d_f,$ 
$\alpha=\frac{d_w}{2},$ $p=2,q=\infty.$ We drop the 
superscript `2' in the definition of $b_m$'s.  The
Lipschitz norm of $f,$ with constant `$L$' in place of `2'
will be denoted by  $\|f\|_{\cal L}$ (as explicitly written below).

After this technical lemma we can pass to our theorem.
\begin{th}
Let $\alpha=\frac{d_w}{2}$ and let 
 $\cll = {Lip}\,(\alpha,2,\infty)(\ck),$ 
endowed with the norm 
\[\|f\|_{\cal L}=\|f\|_2+\|(b_m(f))\|_{\infty}.\]
Let $\cad$
be the domain of the Dirichlet form on $\ck$ with the norm
(\ref{i8}).
Then $\cll = \cad $ and the norms $\|f\|_{\cal L}$
and $\|f\|_{\ce}$ are equivalent.
\end{th}

\subparagraph*{Proof.}
\noindent{\bf Part 1.} First we prove that 
 for $f\in\cad$
we have
$\|f\|_{L}\leq C\|f\|_{\ce}.$ To this end, it is enough to see
that there exists a constant $c>0$ such that for each $f\in{\cal D}$ and
$m=1,2,...$ 
\[\|b_m(f)\|^2_{\infty}\leq c\cdot \ce(f,f).\]
Take $f\in \cad.$
Then
\[\ce(f,f)=\lim_{m\to\infty} c\cdot \left(\frac{L^{d_w}}{L^{d_f}}\right)^m
\sum_{S\in \cf_m}\sum_{x,y\in V(S)} a^{(m)}_{x,y}(f(x)-f(y))^2<+\infty,\]
for some constant $c.$
As all the conductivities $a^{(m)}_{x,y}$ are positive and
 chosen from a finite set, we have
\[\lim_{m\to\infty}\left(\frac{L^{d_w}}{L^{d_f}}\right)^m\sum_{S\in\cf_m}
\sum_{x,y\in V(S)}(f(x)-f(y))^2\leq  c\cdot \ce(f,f).\]

Let $m>0$ be fixed. We prove that
for $c_0=\inf\{|x-y|: x,y\in V_0\},$ for each $m\geq 1$
\begin{equation}\label{eq4}
(b_m(f))^2 =
L^{m(d_w+d_f)}{\int\!\!\int}_{|x-y|<\frac{c_0}{L^m}}(f(x)-f(y))^2
d\mu(x)d\mu(y)\leq c\cdot \ce(f,f).
\end{equation}

Recall that  the normalized counting measures
 on the sets $V_n,$ $\mu_n$'s,
converge weakly to the $x^{d_f}-$Hausdorff measure on the
fractal, $\mu.$ For each $x\in\ck$ and $r>0,$
  $\mu(\partial(B(x,r)))=0.$ 
  Moreover $\cad \subset C(\ck),$ so that $f$ in question is continuous
  (uniformly continuous in fact). As a consequence of the 
  Portmanteau lemma 
  it is enough to see that
   (\ref{eq4}) holds for the approximations
$\mu_n,$ i.e. that
\[
\forall_{n>m}\;\;
L^{m(d_w+d_f)}{\ds\int\!\!\ds\int}_{|x-y|<\frac{c_0}{L^m}}(f(x)-f(y))^2
d\mu_n(x)d\mu_n(y)\leq c\cdot\ce(f,f),
\]
with constant not depending  on $n,m,f.$

For $S\in \cf_m,$ let $S_{*}$ be the union of $S$ and all those symplices
from $\cf_m$ that have a point in common with $S.$ 
Our choice of $c_0$ ensures that if $x\in S\in\cf_m$ and 
$|x-y|<\frac{c_0}{L^m},$ then $y\in S_{*}.$ Thus
\begin{eqnarray}\label{eq5}
{\int\!\!\int}_{|x-y|<\frac{c_0}{L^m}}(f(x)-f(y))^2
 d\mu_n(x)d\mu_n(y)\leq\nonumber\\
\leq\frac{1}{2} \frac{1}{[\#V_n]^2}
\sum_{S\in\cf_m}\sum_{x\in S\cap V_n}\sum_{y\in S_{*}\cap V_n}
(f(x)-f(y))^2.
\end{eqnarray}

There exists 
$z_{xy}\in V(S)$ such that $z_{xy}$ and $x,$ as well as $z_{xy}$ and $y$
 belong to the same $m$-symplex.
As
\begin{equation}\label{eq6}
(f(x)-f(y))^2\leq 2[(f(x)-f(z_{xy}))^2+(f(z_{xy})-f(y))^2],
\end{equation}

(\ref{eq5}) is not greater than
\begin{eqnarray}\label{eq8}
\frac{c}{(\# V_n)^2}\sum_{S\in\cf_m}
\left(\sum_{x\in S\cap V_n}\sum_{y\in S_{*}\cap V_n}[(f(x)-f(z_{xy}))^2+
(f(z_{xy})-f(y))^2]\right)\nonumber\\
\leq \frac{c}{(\# V_n)^2}\sum_{S\in\cf_m}\sum _{x\in S\cap V_n}
\left(\sum_{z\in V(S)} (f(x)-f(z))^2\cdot \# [V_n\cap S]\right)\nonumber\\
=\frac{c\cdot\# V_{n-m}}{(\# V_n)^2}\sum_{S\in\cf_m}\left(\sum_{z\in V(S)}
\sum_{x\in S\cap V_n}(f(x)-f(z))^2\right).
\end{eqnarray}

\noindent
Let $m\geq 1$ and  $S\in \cf_m$  be fixed.
$S$ has unique address $(i_1,i_2,...,i_m)$, i.e.
\[S=\psi_{i_1}\circ\dots\circ\psi_{i_m}(\ck),\;\;\mbox{ for some }
\;i_1,...,i_m\in\{1,...,M\}.\]
Its vertices $z\in V_m$ are of the form 
\[z=\psi_{i_1}\circ\dots\circ\psi_{i_m}(y)\]
and any $x\in S\cap V_n$  ($n>m)$
can be written as
\[x=\psi_{i_1}\circ\dots\circ\psi_{i_m}\circ\psi_{i_{m+1}}
\circ\dots\psi_{i_n}(\bar y),\]
for some $y,\bar y\in V_0.$
It follows that the double sum in round brackets  in (\ref{eq8}) is not greater
than (there is no
equality
as there exist points with multiple addresses)
\begin{equation}\label{eq7}
\sum_{y\in V_0}\sum_{i_{m+1},...,i_{n}}
\sum_{\bar y\in V_0}
((f(\psi_{i_1}\circ\dots\circ\psi_{i_m}(y))-
f(\psi_{i_1}\circ\dots\circ\psi_{i_n}(\bar y)))^2.
\end{equation}
Next, we build a `path' from $z$ to $x$ as follows:
\[\begin{array}{rcl}
s_0 & = & \psi_{i_1}\circ\dots\circ\psi_{i_m}(y)=z\in V_m\cap S,\\
s_1 & = & \psi_{i_1}\circ\dots\circ\psi_{i_m}\circ\psi_{i_{m+1}}(y)\in
 V_{m+1}\cap S,\\
\multicolumn{3}{c}{\vdots}\\
s_{n-m-1} & = & \psi_{i_1}\circ\dots\circ\psi_{i_{n-1}}(y)\in V_{n-1}\cap
S,\\
s_{n-m} &  = & \psi_{i_1}\circ\dots\circ\psi_{i_n}(\bar y)=x\in V_n\cap S.
\end{array}\]
Using (\ref{eq6}) repeatedly we get
\[\begin{array}{l}
(f(z)-f(x))^2)\leq\\
\multicolumn{1}{r}{\leq 2(f(s_0)-f(s_1))^2+4(f(s_1)-f(s_2))^2+...+
2^{n-m}(f(s_{n-m-1})-f(s_{n-m}))^2.}\end{array}\]
Observe that for $l=0,1,...,n\!-m\!-1$ both points $s_l$ and $s_{l+1}$
belong
to the $(m+l)$-symplex $S_l=\psi_{i_1}\circ\dots\circ\psi_{i_{m+l}}(\ck),$
and that $s_l\in V(S_l).$
How many times can a given pair $(s_l,s_{l+1})$ of this form
 appear in such a path? 
It appears whenever the addresses of $x$ and $s_{l+1}$ coincide up to
the $i_{m+l}-$th place, i.e. at most $M^{n-(l+1)}$ times. Summing up
we arrive at:
\[
(\ref{eq7})\leq c\cdot\sum_{r=m}^{n-1}\sum_{ \cf_r\ni  R\subset S}
\sum_{z\in V(R)}\sum_{w\in R\cap V_{r+1}} 2^{r-m+1}M^{n-(r+1)}
(f(z)-f(w))^2.\]
Every symplex $R\in\cf_r$ is just a scaled-down copy of $\ck,$ therefore
there exists a constant depending on $\ck$ only such that
\[\sum_{z\in V(R),w\in R\cap V_{r+1}}(f(z)-f(w))^2\leq 
c\cdot \widetilde{\ce}^{(r+1)}_R(f,f),\]
with
$\widetilde{\ce}^{(r+1)}_R(f,f)$
defined by (\ref{aa}).
To see this, connect $x$ and $z$ 
by a path $(p_0,p_1,...,p_a)$ as follows:
$p_0=x,p_a=z,$  $p_i\in V_{r+1}\cap R$ for $i=0,1,...,a$ and for each
$i=0,...,r-1$ $(p_i,p_{i+1})$ are $(r+1)$-neighbours.
Next
use the inequality $(a_1+...+a_p)^2\leq p(a_1^2+...+a_p^2),$ valid for
arbitrary real numbers $a_i$ and positive integer $p,$ then
observe that a given pair  $(u,v)$ of $(r+1)$-neighbours can appear
in such a path only a finite number of times and that this number does not
change with $r$ ($ R$ is just a scaled-down copy of $\ck$ and
${R}\cap V_{r+1}$ --- a scaled-down copy of $V_1).$

Summing over $\cf_r\ni R\subset S$ and then over $S\in \cf_m$ we end up with
\[\begin{array}{ll}
(\ref{eq8}) & \leq c\cdot\frac{M^{n-m}}{M^{2n}}\ds\sum_{S\in\cf_m}
\ds\sum_{r=m}^{n-1}
\ds\sum_{\cf_r\ni R\subset S} 2^{r-m+1}M^{n-(r+1)}
\widetilde{\ce}^{(r+1)}_R(f,f)\\
{} & =c\cdot\frac{M^{n-m}}{M^{2n}}\ds\sum_{r=m}^{n-1}
2^{r-m+1}M^{n-(r+1)}\ds\sum_{R\in \cf_r}\widetilde{\ce}^{(r+1)}_R
(f,f)\\
{} & = c\cdot\frac{M^{n-m}}{M^{2n}}\ds\sum_{r=m}^{n-1}
2^{r-m+1}M^{n-(r+1)}\widetilde{\ce}^{(r+1)}(f,f).
\end{array}\]
Recall that the sequence 
\[\left(\frac{L^{d_w}}{L^{d_f}}\right)^r
\widetilde{\ce}^{(r)}(f,f)\] increases 
towards $c\cdot\ce(f,f) $ when $r\to\infty,$
 so that
\[(\ref{eq8})\leq c\cdot\ce (f,f)\left[\frac{1}{M^m}\sum_{r=m}^{n-1}
2^{r-m+1}\frac{1}{M^{r+1}}
\left(\frac{L^{d_f}}{L^{d_w}}\right)^{r+1}\right].\]
After summing up the resulting geometric series, using $L^{d_f}=M,$
we obtain that the term in square brackets is not bigger than
\[\frac{1}{M^m}\cdot\frac{1}{L^{md_w}}\cdot\frac{1}{1-\frac{2}{L^{d_w}}}=
\frac{c}{(L^{d_w+d_f})^m},\]
if $L^{d_w}>2.$ 
But as long as $d_s <2$ (equivalently: $\rho >1)$  and $M\geq 2$ we have
$L^{d_w}=
\rho L^{d_f}=\rho M>2 $ (see the discussion at the end of
section \ref{constants}).
Proof of the first part is completed.

\noindent{\bf Part 2. } Now we show that 
$\cll\subset\cad.$
Suppose that
$g\in Lip\,(\alpha,2,\infty)(\ck),$ i.e. $g\in L^2(\ck, \mu),$
and that
$\|(b_m(g))\|_{\infty} <+\infty,$ with $b_m(g)$ defined by (\ref{eq1}).
We need to show that 
\[\sup_{m\geq 1}\left(\frac{L^{d_w}}{L^{d_f}}\right)^m
\widetilde{\ce}^{(m)}(g,g)
<c\cdot\|(b_m(g))\|^2_{\infty}.\]
Again,
\begin{equation}\label{eq2.0}
\begin{array}{ll}
\left(\frac{L^{d_w}}{L^{d_f}}\right)^m\widetilde{\ce}^{(m)}(g,g)  &
 = c\cdot \left(\frac{L^{d_w}}{L^{d_f}}\right)^m\ds\sum_{S\in\cf_m}
\ds\sum_{x,y\in V(S)}a^{(m)}_{x,y}(g(x)-g(y))^2 \\
{} & \leq c\cdot \left(\frac{L^{d_w}}{L^{d_f}}\right)^m
\ds\sum_{S\in\cf_m}\ds\sum_{x,y\in V(S)}(g(x)-g(y))^2.
\end{array}\end{equation}
Fix $m\geq 1$ and $S\in \cf_m.$ Then for each $p\in S$
\begin{equation}\label{l1}
(g(x)-g(y))^2\leq 2[(g(x)-g(p))^2+(g(p)-g(y))^2].
\end{equation} 
Integrating both sides of (\ref{l1})
 over $p\!\in\! S$ with respect to the Hausdorff
measure $\mu$ we get
\[(g(x)-g(y))^2\leq \frac{1}{\mu(S)}\left[\int_S(g(x)-g(p))^2d\mu(p)
+\int_S (g(y)-g(p))^2d\mu(p)\right]\]
so that
\begin{equation}\label{eq2.2}
\sum_{x,y\in V(S)}(g(x)-g(y))^2\leq \frac{2\#[V(S)]}{\mu(S)}\sum_{x\in V(S)}
\int_S(g(x)-g(p))^2d\mu(p).
\end{equation}

Let $(i_1,...i_m)$ be the address of the symplex $S,$ i.e.
$S=\psi_{i_1}\circ...\circ\psi_{i_m}(\ck).$ Let a vertex $x\!\in\! V(S)$ 
be given. It is of the
 form $x=\psi_{i_1}\circ...\circ\psi_{i_m}(v)$ for some $v\in V_0.$
Let $l\in\{1,...M\}$ be such an index that $v=\psi_l(v).$ 
Consider
the following decreasing sequence of symplices:
\[\begin{array}{rcl}
\cf_m & \ni  & S_0=\psi_{i_1}\circ...\circ\psi_{i_m}(\ck)=S,\\
\cf_{m+1} & \ni & S_1=\psi_{i_1}\circ...\circ\psi_{i_m}\circ\psi_l(\ck),\\
\cf_{m+2} & \ni &  S_2=\psi_{i_1}\circ...\circ\psi_{i_m}\circ\psi_l
\circ\psi_l(\ck),\\
\multicolumn{3}{c}{\vdots}\\
\cf_{m+r}& \ni & 
S_r=\psi_{i_1}\circ...\circ\psi_{i_m}\circ(\psi_l)^r(\ck),\\
\multicolumn{3}{c}{\vdots}
\end{array}\]
It is clear that for any $r=0,1,2...,$ $x\in V(S_r)$ 
and that $\bigcap_{r=0}^{\infty}=\{x\}.$

Suppose $k\geq 1$ is given (to be chosen later on). Let
\begin{equation}\label{eq2.1}
T_0=S_0,\,\,T_1=S_k,\,\,T_2=S_{2k},...\;\;\mbox{and so on}.
\end{equation}
For any given $p_i\in T_i,$ $i=1,2,..,\nu$ ($\nu$ --- an arbitrary positive
integer) we have 
\[\begin{array}{c}
(g(x)-g(p_0))^2\leq 2((g(x)-g(p_{\nu}))^2+ (g(p_{\nu})-g(p_0))^2\leq\\[2mm]
\leq 2(g(x)-g(p_{\nu}))^2+[4(g(p_0)-g(p_1))^2+8(g(p_1)-g(p_2))^2+\cdots
+\\[2mm]
+2^{\nu+1}(g(p_{\nu-1})-g(p_{\nu}))^2].
\end{array}\]
Integrating this inequality $(\nu+1)$ times over $p_0\in T_0,$ $p_1\in T_1,$
...$p_{\nu}\in T_{\nu}$ and dividing both sides by $\mu(T_1)\cdot...\cdot
\mu(T_{\nu})$ we obtain:
\[\begin{array}{l}
\ds\int_{T_0}(g(x)-g(p_0))^2d\mu(p_0)\leq
\frac{2\mu(T_0)}{\mu(T_{\nu})}\ds\int_{T_{\nu}}(g(x)-g(p_{\nu}))^2d\mu(p_{\nu})\\
+\ds\sum_{r=0}^{\nu-1}2^{r+2}\frac{\mu(T_0)}{\mu(T_r)\mu(T_{r+1})}
\ds\int_{T_r}\ds\int_{T_{r+1}}(g(p_r)-g(p_{r+1}))^2d\mu(p_r)d\mu(p_{r+1}).
\end{array}\]
As $(T_r)_{r\geq 0}$ is a decreasing sequence of symplices, for each 
$r\geq 0$ both $p_r$ and $p_{r+1}$ belong to $T_r,$ so that 
$|p_r-p_{r+1}|\leq \mbox{diam}\, T_r=\frac{\mbox{\scriptsize diam}\,
{\normalsize\ck}}{L^{m+rk}}.$
Therefore we will
not destroy the upper bound if we replace the integral over $\{p_r\in T_r,
p_{r+1}\in T_{r+1}\}$ by an integral over $\{p_r\in S, p_{r+1}\in \ck,
|p_r-p_{r+1}|<\frac{\mbox{\scriptsize diam}\, {\normalsize\ck}}{L^{m+rk}}\}.$  Also,
as
$\mu(T_r)=
\frac{1}{M^{m+rk}},$ we get
\begin{eqnarray}\label{eq2.3}
\int_{T_0}(g(x)-g(p))^2d\mu(p)\leq 2M^{\nu k}
\ds\int_{T_{\nu}}(g(x)-g(p))^2d\mu(p)+\nonumber\\
+\ds\sum_{r=0}^{\nu-1}2^{r+2}M^{m+(2r+1)k}\ds\int_{S}
\ds\int_{|p-q|<\frac{\mbox{\scriptsize diam}\,
{\small\ck}}{L^{m+rk}}}(g(p)-g(q))^2d\mu(p)d\mu(q).
\end{eqnarray}
Combining (\ref{eq2.2}) and (\ref{eq2.3}) (recall that $T_0=S$ and
 $\mu(S)=\frac{1}{M^m})$ we obtain
\begin{eqnarray}\label{eq2.4}
\sum_{S\in\cf_m}\sum_{x,y\in V(S)}(g(x)-g(y))^2\leq c\cdot \sum_{S\in\cf_m}
\sum_{x\in V(S)}\left( 2M^{\nu k}\int_{T_{\nu}}(g(x)-g(p))^2d\mu(p) \right.
\nonumber\\
\left. +\sum_{r=0}^{\nu-1}2^{r+2}M^{m+(2r+1)k}\!\int_S\int_{|p-q|<\frac
{\mbox{\scriptsize diam}\,
{\small\ck}}{L^{m+rk}}}(g(p)-g(q))^2d\mu(p)d\mu(q)\!\right).
\end{eqnarray}

As $x\in V(T_{\nu}), $ we have $|x-p_{\nu}|\leq \mbox{diam}\,T_{\nu}=\frac
{\mbox{\small diam}\,{\small\ck}}{L^{m+\nu k}}.$ From Theorem \ref{continuity}
\linebreak
$g$ is H\"older continuous with exponent $\frac{d_w-d_f}{2},$
as long as $\frac{d_w-d_f}{2}\in(0,1).$
This is true for all simple nested fractals: $d_w>d_f$ is equivalent to
$d_s<2.$ On the other hand we always have $d_w-d_f\leq 1$
(otherwise the resolvent density of the Brownian motion,
 being H\"older continuous with
exponent $d_w-d_f,$ would be a constant function). 
Therefore, in view of  $L^{d_f}=M,$
\[\int_{T_{\nu}}(g(x)-g(p))^2d\mu(p)\leq c\cdot \mu(T_{\nu})
(\mbox{diam}\,T_{\nu})^{d_w-d_f}=c\cdot
\left(\frac{1}{L^{m+\nu k}}\right)^{d_w}.\]
It follows for each $\nu, k\geq 1$
\begin{equation}\label{eq2.5}
\begin{array}{c}
\ds\sum_{S\in\cf_m}\!\ds\sum_{x,y\in V(S)}(g(x)-g(y))^2\leq
c M^m\ds\sum_{S\in \cf_m}\!\ds\sum_{x\in V(S)} 2 M^{\nu k}\ds\int_{T_{\nu}}
(g(x)-g(p))^2 d\mu(p) + \\[2mm]
+c M^m \!\!\ds\sum_{S\in \cf_m}\!\# [V(S)]\ds\sum_{r=0}^{\nu-1}M^{m+(2r+1)k}
2^{r+2}\!\ds\int_{p\in S}\!\ds\int_{|p-q|<\frac{\mbox{\scriptsize
diam}\,\ck}{L^{m+rk}}}
(g(p)-g(q))^2d\mu(q)d\mu(p).
\end{array}
\end{equation}
Finally, as $\# \cf_m = M^m,$ using (\ref{eq2.0}),  (\ref{eq2.5}) and replacing
$p\in  S$ by $p\in \ck$ in the last double integral above, we have
\[\begin{array}{l}
\left(\frac{L^{d_w}}{L^{d_f}}\right)^m\widetilde{\ce}^{(m)}(g,g)\leq\\
\leq c\cdot \frac{L^{m d_w}}{M^m}\left(
\#\cf_m\cdot M^m\left(\frac{1}{L^{m+\nu
k}}\right)^{d_w}+M^m\ds\sum_{r=0}^{\nu-1}
M^{(r+1)k}M^{m+rk}2^{r+2}\frac{(b_{m+rk}(g))^2}
{L^{(m+rk)(d_w+d_f)}}\right)\leq\\
\leq c\cdot \left(\frac{M^m}{L^{\nu k d_w}}+\|(b_m(g))\|^2_{\infty}
\ds\sum_{r=0}^{\nu-1}\left(\frac{2M^k}{L^{kd_w}}\right)^r\right).
\end{array}\]

As $\frac{M^k}{L^{kd_w}}=(L^{d_w-d_f})^k$ and $d_w>d_f,$ we can choose 
$k$ so big that $\left(\frac{2M^k}{L^{kd_w}}\right) <\frac{1}{2}$
and the resulting geometric series is convergent.
 For this choice of $k$ we have,
for all $m,\nu\geq 1:$
\[\left(\frac{L^{d_w}}{L^{d_f}}\right)^m\widetilde{\ce}^{(m)}(g,g)\leq
c\cdot \left(\frac{M^m}{L^{\nu k d_w}}+\|(b_m(g))\|^2_{\infty}\right).
\]
The constant appearing in this estimate depends on the fractal only.
 Therefore,
for fixed $m,$ we can pass to infinity with $\nu,$ getting that
\[\left(\frac{L^{d_w}}{L^{d_f}}\right)^m\widetilde{\ce}^{(m)}(g,g)\leq
c\cdot \|(b_m(g))\|^2_{\infty},\]
and finally
\[\ce(g,g)\leq c\cdot \|(b_m(g))\|^2_{\infty}.\]
The proof is complete.
$\Box$

\end{document}